\documentclass[a4paper,12pt,leqno]{article}

\usepackage{amsmath,amsfonts,latexsym}
\usepackage{amssymb}
\usepackage{graphics}

\begin{document}
%\begin{article}

\renewcommand{\theequation}{\thesection.\arabic{equation}}

\newcommand{\nc}{\newcommand}
\nc{\tr}{{\triangle}} \nc{\vth}{{\vartheta}}
 \nc{\bt}{{\beta}}
\nc{\dl}{{\delta}} \nc{\Dl}{{\Delta}}
 \nc{\p}{{\psi}}
\nc{\gm}{{\gamma}} \nc{\Gm}{{\Gamma}} \nc{\sg}{{\sigma}}
\nc{\ve}{{\varepsilon}} \nc{\ch}{{\cal H}} \nc{\cf}{{\cal F}}
\nc{\cp}{{\cal P}}
 \nc{\td}{\tilde}

\newtheorem{rem}{Remark}[section]
\newtheorem{Th}{Theorem}[section]
\newtheorem{cor}{Corollary}[section]
\newtheorem{Lemma}{Lemma}[section]

%\begin{opening}
\title{Recursive Parameter Estimation: Convergence}
% \subtitle{Basic Instructions}
 \author{Teo Sharia}

\date{}
\maketitle

\begin{center}
{\it
Department of Mathematics \\Royal Holloway,  University of London\\
Egham, Surrey TW20 0EX \\ e-mail: t.sharia@rhul.ac.uk }
\end{center}
%\date: rather not

%\dedication{To Jim}

%\translation{De Kluwer LaTeX stylefile; aanwijzingen voor auteurs}

%\runningtitle{Recursive Parameter Estimation}
%\runningauthor{T. Sharia}

%\begin{ao}
%Department of Mathematics \\Royal Holloway,  University of London\\
%Egham, Surrey TW20 0EX \\ e-mail: t.sharia@rhul.ac.uk
%\end{ao}

%\begin{motto}
%What can't be done with TeX isn't worth doing.
%\end{motto}

\begin{abstract}
We consider estimation procedures which are recursive in the sense
that each successive
 estimator is obtained from the previous one by a simple adjustment. We propose
  a wide class of recursive estimation
       procedures for the general statistical model and study convergence.
\end{abstract}

\begin{center}
\end{center}

\begin{center}
%\end{opening}
\end{center}

%***************************************************************

%***************************************************************************

                       %S E C T I O N 1
\section{Introduction}
Let $X_1, \dots, X_n$ be independent identically distributed
(i.i.d.) random variables (r.v.'s) with  a common distribution
function
  $F_{\theta}$ with a  real unknown parameter $\theta$. An
$M$-estimator of $\theta$ is defined as a statistic $
\hat\theta_n=\hat\theta_n(X_1,\dots,X_n), $ which is a solution
w.r.t. $v$ of the estimating equation
%%%%%%%%%%%%%%%%%%%%%%%%%%%%%%%%%%%%%%%%%%%%%%%%%
%                       (esteq)
%%%%%%%%%%%%%%%%%%%%%%%%%%%%%%%%%%%%%%%%%%%%%%
\begin{equation} \label{esteq}
\sum_{i=1}^n \psi(X_i;v)=0,
\end{equation}
where $\psi$ is a suitably chosen function. For example, if
$\theta$ is a location parameter in  the normal family of
distribution functions, the choice   $\psi(x,v)=x-v$ gives
 the MLE (maximum likelihood estimator).
For the same problem, if $\psi(x,v)=\mbox{sign}(x-v),$ the
solution of \eqref{esteq} reduces to the median of $X_1, \dots,
X_n$. In general, if $f(x,\theta)$ is the probability density
function (or probability function) of  $F_{\theta}(x)$  (w.r.t. a
$\sigma$-finite measure $\mu$)  then the choice
$\psi(x,v)=f'(x,v)/f(x,v)$ yields the MLE.

 Suppose now that $X_1, \dots, X_n$ are not
necessarily independent or identically  distributed r.v's, with a
joint distribution depending on a real parameter $\theta$. Then an
$M$-estimator of $\theta$ is defined as  a solution of the
estimating equation
%%%%%%%%%%%%%%%%%%%%%%%%%%%%%%%%%%%%%%%%%%%%%%%%%%%%
%                  \eqref{esteqg}
%%%%%%%%%%%%%%%%%%%%%%%%%%%%%%%%%%%%%%%%%%%%%%%%%%%%%%
\begin{equation}\label{esteqg}
\sum_{i=1}^n \psi_i(v)=0,
\end{equation}
where   $\psi_i(v)=\psi_i(X_{i-k}^i; v)$ with $X_{i-k}^i=
(X_{i-k}, \dots,X_i )$. So, the $\psi$-functions  may now depend
on the past observations as well. For instance, if $X_i$'s are
observations from a discrete time Markov process, then one can
assume that   $k=1$. In general, if no restrictions are placed on
the dependence structure of the process $X_i$, one may need to
consider $\psi$-functions depending on the vector of all past and
present observations of the process (that is, $k=i-1$). If the
conditional probability density function (or probability function)
of the observation $X_i,$ given $X_{i-k}, \dots,X_{i-1},$ is
$f_i(x,\theta)=f_i(x,\theta|X_{i-k}, \dots,X_{i-1})$,  then one
can obtain the  MLE  on choosing
$\psi_i(v)=f'_i(X_i,v)/f_i(X_i,v).$ Besides MLEs, the class of
$M$-estimators includes estimators
 with special properties such as {robustness}.
 %An estimator
 %is  robust if its behaviour is not ``seriously
 %affected'' by violations of underlying assumptions.
Under certain regularity and ergodicity conditions it can be
proved that there exists a consistent sequence of solutions of
\eqref{esteqg} which has  the property of local asymptotic
linearity. (See e.g., Serfling \cite{Ser},
Huber \cite{Hub}, Lehman \cite{Leh}.  A comprehensive bibliography can be
found in Launer and Wilkinson \cite{Laun}, Hampel at al \cite{Ham}, Rieder
\cite{Ried}, and Jure$\check{c}$kov$\acute{a}$ and Sen \cite{Jur}.)

 If   $\psi$-functions are
 nonlinear, it is rather difficult to
 work with  the corresponding estimating equations, especially
  if for every sample size $n$
 (when new data are acquired),
 an estimator has to be computed afresh. In this paper we consider
 estimation procedures which are recursive in the sense that each successive
 estimator is obtained from the previous one by a simple adjustment.
 Note that for a linear estimator,
 e.g., for the sample mean, $\hat\theta_n=\bar X_n$ we have
$\bar X_n=(n-1)\bar X_{n-1}/n+X_n/n$, that is
$\hat\theta_n=\hat\theta_{n-1}(n-1)/n+X_n/n$, indicating that the
estimator $\hat\theta_n$ at  each  step $n$ can be obtained
recursively using the estimator at the previous step
$\hat\theta_{n-1}$ and the new information $X_n$. Such an exact
recursive relation may not hold for nonlinear estimators (see,
e.g., the case of the median).

In general, the following heuristic argument can be used to
establish   a possible form of an approximate recursive relation
(see also Jure$\check{c}$kov$\acute{a}$ and Sen \cite{Jur},
 Khas'minskii  and Nevelson \cite{Khas},
 Lazrieva and  Toronjadze \cite{Laz3}).
 Since  $\hat\theta_n$ is defined
as   a root  of the estimating equation \eqref{esteqg}, denoting
the left hand side of \eqref{esteqg} by $M_n(v)$ we have
$M_n(\hat\theta_n)=0$ and  $M_{n-1}(\hat\theta_{n-1})=0$. Assuming that
the difference $\hat\theta_n-\hat\theta_{n-1}$ is ``small'' we can
write
$$
0=M_n(\hat\theta_n)-M_{n-1}(\hat\theta_{n-1})=
M_n\left(\hat\theta_{n-1}+(\hat\theta_n-\hat\theta_{n-1})\right)-
M_{n-1}(\hat\theta_{n-1})
$$
$$
\approx
M_n(\hat\theta_{n-1})+M_n'(\hat\theta_{n-1})(\hat\theta_n-\hat\theta_{n-1})
-M_{n-1}(\hat\theta_{n-1})
$$
$$
=M_n'(\hat\theta_{n-1})(\hat\theta_n-\hat\theta_{n-1})
+\psi_n(\hat\theta_{n-1}).
$$
Therefore,
$$
\hat\theta_n\approx\hat\theta_{n-1}-\frac{\psi_n(\hat\theta_{n-1})}
{M_n'(\hat\theta_{n-1})},
$$
where $M_n'(\theta)=\sum_{i=1}^n \psi'_i(\theta)$. Now,    depending
 on the nature of the underlying model,  $M_n'(\theta)$ can be replaced by
 a simpler expression. For instance, in  i.i.d.
 models with $\psi(x,v)=f'(x,v)/f(x,v)$  (the  MLE case), by the strong law of large
 numbers,
 $$
 \frac {M_n'(\theta)}n =\frac1n  \sum_{i=1}^n
 \left(f'(X_i,\theta)/f(X_i,\theta)\right)'\approx E_\theta
 \left[  \left(f'(X_1,\theta)/f(X_1,\theta)\right)'\right]
 =-i(\theta)
 $$
 for large $n$'s,
where $i(\theta)$ is the one-step Fisher   information. So, in
this case, one can use the recursion

%%%%%%%%%%%%%%%%%%%%%%%%%%%%%%%%%%%%%%%%%%%%%%%%%%
%                 {mliid}
%%%%%%%%%%%%%%%%%%%%%%%%%%%%%%%%%%%%%%%%%%%%%%%%%
\begin{equation} \label{mliid}
\hat\theta_n=\hat\theta_{n-1}+\frac{1}{n \; i(\hat\theta_{n-1})}
\frac{f'(X_n,\hat\theta_{n-1})}{f(X_n,\hat\theta_{n-1})},
~~~~~~~~~n\ge 1,
\end{equation}
to construct an estimator which is ``asymptotically equivalent''
to the MLE.

Motivated by  the above argument,   we consider  a class of
estimators
%%%%%%%%%%%%%%%%%%%%%%%%%%%%%%%%%%%%%%%%%%%%%%%%%%
%                   \eqref{rec1}
%%%%%%%%%%%%%%%%%%%%%%%%%%%%%%%%%%%%%%%%%%%%%%%%%
\begin{equation} \label{rec1}
\hat\theta_n=\hat\theta_{n-1}+
{\Gamma_n^{-1}(\hat\theta_{n-1})}\psi_n(\hat\theta_{n-1}),
~~~~~~~~~n\ge 1,
\end{equation}
where $\psi_n$ is a suitably chosen vector process, $\Gamma_n$  is
a (possibly random) normalizing matrix process and
$\hat\theta_0\in {\mathbb{R}}^m$ is some initial value. Note that
while the main goal is  to study recursive procedures with
non-linear $\psi_n$ functions, it is worth mentioning that any
linear estimator can be written in  the form  \eqref{rec1} with
linear, w.r.t. $\theta$, $\psi_n$ functions. Indeed, if
$\hat\theta_n=\Gamma_n^{-1}\sum_{k=1}^n h_{k} (X_k),$  where
$\Gamma_k$ and $h_k(X_k)$ are  matrix and vector processes  of
suitable dimensions, then (see Section 4.2 for details)
$$
\hat\theta_n=\hat\theta_{n-1}+\Gamma_n^{-1}\left( h_n(X_n)-
(\Gamma_n-\Gamma_{n-1})\hat\theta_{n-1}\right),
$$
which is obviously of the form  \eqref{rec1} with $\psi_n(\theta)=
h_n(X_n)- (\Gamma_n-\Gamma_{n-1})\theta.$

It should be noted that at first glance, recursions \eqref{mliid} and  \eqref{rec1}
resemble the  Newton-Raphson iterative procedure  of  numerical optimisation.
In the i.i.d. case,  the Newton-Raphson iteration for the likelihood
equation  is
%%%%%%%%%%%%%%%%%%%%%%%%%%%%%%%%%%%%%%%%%%%%%%%%%
%                 {NR}
%%%%%%%%%%%%%%%%%%%%%%%%%%%%%%%%%%%%%%%%%%%%%%%%%
\begin{equation} \label{Nr}
\vartheta_k=\vartheta_{k-1}+J^{-1}(\vartheta_{k-1})
\sum_{i=1}^n\frac{f'(X_i,\vartheta_{k-1})}{f(X_i,\vartheta_{k-1})},
~~~~~~~~~k\ge 1,
\end{equation}
where $J(v)$ is minus  the second logarithmic derivative of the log-likelihood function, that is,
$-\sum_{i=1}^n\frac{\partial}{\partial v}\left(f'(X_i,v)/f(X_i,v)\right)$
or its expectation, that is, the information matrix  $n i(v).$
In the latter case, the iterative scheme is often called the method of
scoring, see e.g., Harvey \cite{Har}.
(We do not consider  the  so called one-step Newton-Raphson method  since it requires an auxiliary
consistent estimator).
The main feature of the scheme \eqref{Nr} is that $\vartheta_k$,
at each step $k=1,2,\dots,$
is $\sigma (X_1,\dots ,X_n)$ - measurable (where $\sigma (X_1,\dots ,X_n)$ is the $\sigma$-field
generated by the random variables $ X_1,\dots ,X_n$). In other words, \eqref{Nr} is a deterministic
procedure to find a root, say $\tilde\theta_n$, of the likelihood equation\\
$\sum_{i=1}^n \left(f'(X_i,v)/f(X_i,v)\right) =0.$ On the other hand
the random variable $\hat \theta_n$  derived from
\eqref{mliid} is  an estimator of $\theta$ for each n=1,2,\dots (is $\sigma
(X_1,\dots ,X_n)$- measurable at each $n$). Note also that in the iid case,
 \eqref{mliid} can be regarded as a stochastic iterative scheme,
i.e.,  a classical stochastic approximation procedure,  to detect the root of
an unknown function when the latter can only be observed with random errors (see Remark 3.1).
A theoretical implication of this is that by studying the procedures \eqref{mliid},
or in general \eqref{rec1},
we study  asymptotic behaviour of the estimator of the unknown parameter.
As far as applications  are concerned, there are several advantages in  using
\eqref{rec1}. Firstly, these procedures are easy to use since
each successive
 estimator is obtained from the previous one by a simple adjustment
and without  storing all the data unnecessarily. This is especially
convenient when the data come sequentially.
 Another potential benefit of using \eqref{rec1}  is that it allows one to
 monitor and detect certain changes  in probabilistic characteristics  of
 the underlying process such as change of the value of the unknown parameter. So, there
 may be a benefit in using these procedures in linear cases as well.

  In  i.i.d. models,
 estimating procedures similar to \eqref{rec1} have been studied by a number of
authors using methods of stochastic approximation theory (see,
e.g., Khas'minskii  and Nevelson \cite{Khas}, Fabian \cite{Fab},
 Ljung and Soderstrom
\cite{Ljun2}, Ljung et al \cite{Ljun1},  and references therein).
Some work has been done for  non i.i.d. models as well. In
particular,  Englund et al \cite{Eng} give an
asymptotic representation results for  certain type of  $X_n$
processes. In Sharia \cite{Shar1} theoretical results on convergence,
rate of convergence and the asymptotic representation are given
under certain regularity and ergodicity assumptions on the model,
in the one-dimensional  case with
$\psi_n(x,\theta)=\frac{\partial}{\partial\theta} \mbox{log}
f_n(x,\theta)$ (see also Campbell \cite{Cam},
%Sharia (1992),
 Sharia \cite{Shar2} and  Lazrieva et al \cite{Laz1}).

In the present paper, we study multidimensional  estimation
procedures of type \eqref{rec1} for the general statistical model.
Section 2 introduces the basic model, objects  and notation.  In
Section 3, imposing ``global'' restrictions on the processes
$\psi$ and $\Gamma$, we study   ``global'' convergence of the
recursive estimators, that is the convergence  for an arbitrary
starting point $\hat\theta_0$.  In Section 4, we  demonstrate the use of these
results on some examples. (Results on rate of convergence,
asymptotic linearity and efficiency, and numerical simulations   will appear
in subsequent publications, see Sharia \cite{Shar3}, \cite{Shar4}.)

%\newpage

                       %S E C T I O N 2

\section {Basic  model,  notation and preliminaries}
\setcounter{equation}{0}

Let   $X_t,\;\; t=1,2,\dots ,$
 be  observations taking values in a  measurable space
$({\bf X},{\cal B}({\bf X}))$  equipped with  a $\sigma$-finite
measure $\mu.$ Suppose  that the distribution of the process $X_t$
depends on an unknown parameter $\theta \in \Theta,$ where
$\Theta$ is an open subset of the $m$-dimensional Euclidean space
$\mathbb{R}^m$. Suppose also  that for each  $t=1,2,\dots$, there
exists a  regular conditional probability density of $X_t$ given values of
past observations  of $X_{t-1},\dots , X_2, X_1$, which will be
denoted by
$$
f_ t(\theta, x_t \mid x_1^{t-1})= f_ t(\theta, x_t \mid
x_{t-1},\dots ,x_1),
$$
where $ f_ 1(\theta, x_1 \mid x_1^0)= f_ 1(\theta, x_1) $ is the
probability  density  of the random variable $X_1.$
Without loss of generality we assume
that all random variables are defined on a probability space $
 (\Omega, \cf)
 $
and denote by $ \left\{P^\theta, \; \theta\in \Theta\right\} $ the
family of the corresponding distributions on $
 (\Omega, \cf).
 $

  Let  $\cf_t=\sigma (X_1,\dots ,X_t)$ be the $\sigma$-field
generated by the random variables $ X_1,\dots ,X_t.$ By
$\left(\mathbb{R}^m, {\cal B}( \mathbb{R}^m ) \right)$ we denote
the $m$-dimensional Euclidean space with the Borel
$\sigma$-algebra ${\cal B}( \mathbb{R}^m )$. Transposition of
matrices and vectors is denoted by $T$. By $(u,v)$ we denote the
standard scalar product of $u,v \in \mathbb{R}^m,$ that is,
$(u,v)=u^Tv.$

Suppose that $h$ is a real valued function  defined on
 $\Theta\subset  {{\mathbb{R}}}^m$.  We denote by $\dot h(\theta)$  the row-vector
 of partial derivatives
of $h(\theta)$ with respect to the components of $\theta$, that
is,
 $$
 \dot h(\theta)=\left(\frac{{\partial}}{{\partial} \theta^1} h(\theta), \dots,
 \frac{{\partial}}{{\partial} \theta^m} h(\theta)\right).
 $$
Also we denote by  \"{h}$(\theta)$ the  matrix of second partial
derivatives.
 The $m\times m$ identity matrix is denoted by ${{\bf 1}}$.

 If for each  $t=1,2,\dots$, the derivative
$ \dot f_t(\theta, x_t \mid x_1^{t-1})$ w.r.t. $\theta$
  exists, then we
can define the function
$$
l_t(\theta, x_t \mid x_1^{t-1})=\frac 1 {f_ t(\theta, x_t \mid
x_1^{t-1})} \dot f_t^T(\theta, x_t \mid x_1^{t-1})
$$
with the convention $0/0=0$.

The {\it one step   conditional  Fisher information matrix} for
$t=1,2,\dots$ is defined as
$$
i_t(\theta\mid x_1^{t-1}) =\int l_t(\theta,z\mid x_1^{t-1})
l^{T}_t(\theta,z\mid x_1^{t-1})f_t(\theta,z\mid x_1^{t-1})\mu
(dz).
$$
We shall use the  notation
$$
f_t(\theta)=f_t(\theta, X_t \mid X_1^{t-1}), \;\;\; \;\;\;
l_t(\theta)=l_t(\theta, X_t\mid X_1^{t-1}),
$$
$$
i_t(\theta)=i_t(\theta\mid X_1^{t-1}).
$$
Note that  the process $i_t(\theta)$  is {\it ``predictable''},
that is, the random variable $i_t(\theta),$ is $\cf_{t-1}$
measurable for each $t\ge 1.$

 Note also that by   definition,
$i_t(\theta)$ is   a version of the conditional expectation w.r.t.
${\cal{F}}_{t-1},$  that is,
$$
i_t(\theta)= E_\theta\left\{l_t(\theta) l^{T}_t(\theta) \mid
{\cal{F}}_{t-1}\right\}.
$$
Everywhere in the present work conditional expectations  are meant
to be
 calculated as integrals w.r.t. the conditional probability densities.

The {\it conditional  Fisher information} at time $t$ is
$$
I_t(\theta)=\sum_{s=1}^t i_s(\theta),  \;\;\;\;\;\;\;\;\;
t=1,2,\dots.
$$
If the $X_t$'s are independent random variables,  $I_t(\theta)$
reduces to the standard Fisher information matrix.  Sometimes
$I_t(\theta)$ is referred as the incremental expected Fisher
information. Detailed discussion of this concept and related work
appears in
 Barndorff-Nielsen and  Sorensen \cite{Bar}, Prakasa-Rao  \cite{Prak}  Ch.3,
 and Hall and Heyde \cite{Hal}.

We say that  ${\bf \psi}= \{\psi_t(\theta, x_t, x_{t-1}, \dots,
x_1)\}_{t\ge 1}$ is a sequence of estimating functions and write
$\bf \psi \in \Psi$, if for each ${t\ge 1},$ $ \psi_t(\theta, x_t,
x_{t-1}, \dots, x_1)  :
 \Theta \times {\bf X}^t \;\;\to \;\;{\mathbb{R}}^m
$ is a   Borel function.

Let $\psi \in {\bf\Psi}$ and denote $\psi_t(\theta)=\psi_t(\theta,
X_t, X_{t-1},\dots, X_1).$ We write $\bf \psi \in \Psi^M$ if  ~
$\psi_t(\theta)$ is a martingale-difference process for each
$\theta \in \Theta,$ ~ i.e.,~ if $
E_\theta\left\{\psi_t(\theta)\mid {\cal{F}}_{t-1}\right\}=0$ for
each $t=1,2,\dots $ (we assume that the conditional expectations
above are well-defined and ${\cal{F}}_{0}$ is the trivial
$\sigma$-algebra).

Note   that if differentiation of the equation
$$
1=\int  f_t(\theta, z \mid x_1^{t-1}) \mu(dz)
$$
is allowed under the integral sign, then $ \{l_t(\theta, X_t \mid
X_1^{t-1})\}_{t\ge 1}\in {\bf\Psi^M}$.

\vskip+0.5cm

 {\bf Convention} {\it  Everywhere in the present work
 $\theta\in \mathbb{R}^m $ is an arbitrary but fixed value
of the parameter.  Convergence and all relations between random
variables are meant with probability one w.r.t. the measure
$P^\theta$ unless specified otherwise. A sequence of random
variables $(\xi_t)_{t\ge1}$ has some property eventually if for
every $\omega$ in a set $\Omega^\theta$ of $P^\theta$ probability
1,
 $\xi_t$  has this property for all $t$ greater than some
$t_0(\omega)<\infty$.}

%\newpage

  %           3   Main results

\section{Main results}
\setcounter{equation}{0}

Suppose   that $\bf \psi \in \Psi$ and    $\Gm_t(\theta)$,  for
each $\theta\in {\mathbb{R}}^m$, is  a predictable $m\times m$ matrix process
with $ \mbox{det} ~\Gm_t(\theta)\neq0$, $t\ge 1$.
 Consider  the estimator $\hat \theta_t$  defined by
                       %(rec2)
\begin{equation}\label{rec2}
\hat \theta_t=\hat \theta_{t-1}+\Gm_t^{-1}(\hat
\theta_{t-1})\p_t(\hat \theta_{t-1}), \qquad  t\ge 1,
\end{equation}
where $\hat \theta_0 \in {\mathbb{R}}^m$
 is an arbitrary  initial point.

Let $\theta \in {\mathbb{R}}^m $ be an arbitrary but fixed value of
the parameter  and  for any  $u\in {\mathbb{R}}^m$ define
 $$
b_t(\theta,u)=E_\theta \left\{\psi_t(\theta+u)\mid
{\cf}_{t-1}\right\}, \qquad
 t\ge 1.
$$
%=E_\theta \left\{\p_t(\theta+u)-\p_t(\theta)\mid{\cf}_{t-1}\right\}.

%%%%%%%%%%%%%%%%%%%%%%%%%%%%%%%%%%%%%%%%%%%%%%%%%%%%%%%%%%%%%%%%%%%%%%%%%%
                 %T H E O R E M  3.1
%%%%%%%%%%%%%%%%%%%%%%%%%%%%%%%%%%%%%%%%%%%%%%%%%%%%%%%%%%%%%%%%%%%%%%%%%%
\begin{Th}
Suppose that
\begin{description}
\item[(C1)]
 $u^{T}\Gm_t^{-1}(\theta+u)b_t(\theta,u) <  0
 \qquad \mbox{for each} ~~ u\ne 0,  \qquad
 P^\theta\mbox{-a.s.}\footnote{Note that the set of $P^\theta$ probability
 $0$ where the inequalities in (C1) and (C3) are not valid should not depend on $u.$}; $
\item[(C2)] $\;$ for each $\ve\in (0, 1),$
$$ \sum_{t=1}^\infty \inf_{\ve \le \|u\| \le {1/\ve}}
|u^{T}\Gm_t^{-1}(\theta+u)b_t(\theta,u)|
 = \infty, \qquad
 P^\theta\mbox{-a.s.};
$$
\item[(C3)] there exists a  predictable scalar process
$(B_t^\theta)_{t\ge1}$ such that
$$
E_\theta \left\{\|\Gm_t^{-1}(\theta+u)\p_t(\theta+u)\|^2\mid
{\cf}_{t-1}\right\}
 \le  B_t^\theta (1+\|u\|^2)
 $$
for each $ u\in {\mathbb{R}}^m, $ $P^\theta$-a.s., and
 $$
 \sum_{t=1}^\infty B_t^\theta < \infty,  \qquad
 P^\theta\mbox{-a.s.}.
$$
 \end{description}
 Then
 $\hat\theta_t$ is  strongly consistent (i.e.,
$\hat \theta_t \to \theta \;\; P^\theta$-a.s.) for any  initial
value $\hat \theta_0$ .
\end{Th}
%%%%%%%%%%%%%%%%%%%%%%%%%%%%%%%%%%%%%%%%%%%%%%%%%%%%%%%%%%%%%%%%%%%%%%%%%%
We will derive this theorem from  a more general result
 (see the end of the section). Let us first comment  on
the  conditions used here.

%%%%%%%%%%%%%%%%%%%%%%%%%%%%%%%%%%%%%%%%%%%%%%%%%%%%%%%%%%%%%%%%%%%%%%%%%%
                           % R E M A R K 3.1
%%%%%%%%%%%%%%%%%%%%%%%%%%%%%%%%%%%%%%%%%%%%%%%%%%%%%%%%%%%%%%%%%%%%%%%%%%
\begin{rem}
{\rm
 Conditions  (C1), (C2),
and (C3) are natural analogues  of the corresponding assumptions
in theory of stochastic approximation. Indeed, let us consider
the i.i.d. case with
$$
f_ t(\theta, z \mid x_1^{t-1})= f(\theta, z), \;\; \;\;
\p_t(\theta)=\p(\theta, z )|_{z=X_t},
$$
where  $\int \p(\theta,z)f(\theta,z)\mu (dz)=0 $ and $
 \Gm_t(\theta)=t\gm(\theta)
$ for some  invertible  non-random matrix  $\gm(\theta)$. Then
$$
   b_t(\theta, u)=b(\theta, u)=\int
\p(\theta+u,z)f(\theta,z)\mu(\,dz),
$$
implying that $b(\theta, 0)=0$. Denote $ \Dl_t=\hat
\theta_t-\theta $ and rewrite \eqref{rec2} in the form
                         % robmonro
\begin{equation}\label{robmonro}
\Dl_t = \Dl_{t-1}  +  \frac1 t\left(\gm^{-1}(\theta+\Dl_{t-1})
b(\theta,\Dl_{t-1})+\ve_t^\theta\right),
\end{equation}
where
$$
\ve_t^\theta=\gm^{-1}(\theta+\Dl_{t-1}) \left\{
\p(\theta+\Dl_{t-1},X_t)-b(\theta,\Dl_{t-1})\right\}.
$$
Equation \eqref{robmonro} defines  a Robbins-Monro stochastic
approximation procedure that converges to  the solution of the
equation
$$
R^{\theta}(u):= \gm^{-1}(\theta+u)b(\theta,u)=0,
$$ when the values of the function
$R^{\theta}(u)$ can only be observed with zero expectation errors
$ \ve_t^\theta$. Note that in  general,  recursion \eqref{rec2}
cannot be considered in the framework of  classical stochastic
approximation theory (see  Lazrieva et al
    \cite{Laz1}, \cite{Laz2} for  the generalized   Robbins-Monro
     stochastic approximations procedures).
 For the i.i.d. case, conditions (C1), (C2) and (C3)
 can be written as   {\bf (I)} and {\bf (II)} in Section 4,
 which  are  standard  assumptions for  stochastic approximation
procedures of  type \eqref{robmonro} (see, e.g., Robbins and Monro
\cite{Rob1}, Gladyshev \cite{Glad},
 Khas'minskii and Nevelson \cite{Khas}, Ljung and Soderstrom \cite{Ljun2},
    Ljung et al \cite{Ljun1}).
 }
 \end{rem}
%\vskip+0.2cm
%%%%%%%%%%%%%%%%%%%%%%%%%%%%%%%%%%%%%%%%%%%%%%%%%%%%%%%%%%%%%%%%%%%%%%%%%%
                  % R E M A R K 3.2
%%%%%%%%%%%%%%%%%%%%%%%%%%%%%%%%%%%%%%%%%%%%%%%%%%%%%%%%%%%%%%%%%%%%%%%%%%
\begin{rem}
{\rm
To understand how the procedure works, consider the one-dimensional case, denote $\Dl_t=\hat\theta_t-\theta$
 and rewrite \eqref{rec2}
in the form
$$
\Dl_t=\Dl_{t-1}+ \Gm_t^{-1}(\theta+\Dl_{t-1})
\p_t(\theta+\Dl_{t-1}).
$$
Then,
$$
E_\theta \left\{\hat\theta_t-\hat\theta_{t-1}\mid{\cf}_{t-1}\right\}
=E_\theta \left\{\Dl_t-\Dl_{t-1}\mid{\cf}_{t-1}\right\}
= \Gm_t^{-1}(\theta+\Dl_{t-1})b_t(\theta,\Dl_{t-1}).
$$
Suppose now that at time ~$t-1,$ ~ $\hat\theta_{t-1}<\theta,$
that is, $\Dl_{t-1}<0.$
 Then, by (C1),
$\Gm_t^{-1}(\theta+\Dl_{t-1})b_t(\theta,\Dl_{t-1})>0$
implying that
$E_\theta \left\{\hat\theta_t-\hat\theta_{t-1}\mid{\cf}_{t-1}\right\}>0.$
So, the next step $\hat\theta_t$ will be in the direction of $\theta$.
 If  at time ~$t-1,$ ~ $\hat\theta_{t-1}>\theta,$
by the same reason, $E_\theta \left\{\hat\theta_t-\hat\theta_{t-1}\mid{\cf}_{t-1}\right\}<0.$
So, the condition (C1) ensures that, on average, at each step the procedure moves towards
$\theta$. However, the magnitude of the jumps
$\hat\theta_t-\hat\theta_{t-1}$ should decrease, for otherwise, $\hat\theta_t$ may oscillate around
$\theta$ without approaching it. This is guaranteed by (C3).
On the other hand, (C2) ensures that  the jumps do not decrease too
rapidly  to avoid  failure of $\hat\theta_t$ to reach $\theta.$

\vskip+0.2cm
\noindent
Now, let us consider a maximum
likelihood type recursive estimator
$$
\hat \theta_t=\hat \theta_{t-1}+I_t^{-1}(\hat
\theta_{t-1})l_t(\hat \theta_{t-1}), \qquad  t\ge1,
$$
where $l_t(\theta) =\dot f_t^T(\theta, X_t \mid X_1^{t-1})/{f_
t(\theta, X_t \mid X_1^{t-1})}$
  and $
I_t(\theta) $ is the conditional  Fisher information with $\det
I_t(\theta) \ne 0$ (see also \eqref{mliid} for the i.i.d. case).
%As we will see in Section 5, $ I_t(\theta) $ is a ``suitable''
%normalizing sequence for the  recursion with the influence process
%$l_t(\theta).$
%\noindent
By Theorem 3.1,  $ \hat \theta_t$ is strongly consistent if
conditions (C1), (C2)  and (C3) are satisfied  with $l_t(\theta)$
and  $I_t(\theta)$ replacing $\p_t(\theta)$ and $\Gm_t(\theta)$
respectively. On the other hand, if  e.g., in the one-dimensional
case, $ b_t(\theta,u) $ is differentiable at $u=0$ and the
differentiation is allowed under the integral sign, then
$$
\frac{\partial}{\partial u} b_t(\theta, u)\overline{}\mid_{u=0}=E_\theta \left\{\dot l_t(\theta)\mid
{\cf}_{t-1}\right\}.
$$
So, if the differentiation w.r.t. $\theta$ ~ of $E_\theta
\left\{l_t(\theta)\mid {\cf}_{t-1}\right\}=0$ is allowed
under the integral sign,
 $ \frac{\partial}{\partial u} b_t(\theta, u)\mid_{u=0}=-i_t(\theta) $
%(see \eqref{long} with $\psi=l$)
 implying that (C1) always holds for small values of
$u\not=0.$

\vskip+0.2cm
\noindent
Condition (C2) in the i.i.d. case  is   a requirement that the
function $ \gm^{-1}(\theta+u) b(\theta,u)$ is separated
from zero on each finite  interval that does not contain $0$. For
the i.i.d. case with continuous w.r.t $u$  functions $
 b(\theta,u)$ and $i(\theta+u),
$ condition (C2) is an easy consequence of (C1).

\vskip+0.2cm
\noindent
Condition  (C3)  is  a boundedness type  assumption which
restricts the growth  of $ \p_t(\theta) $
 w.r.t. $\theta$ with certain
  uniformity w.r.t. $t$.
  }
  \end{rem}
\vskip+0.2cm
%%%%%%%%%%%%%%%%%%%%%%%%%%%%%%%%%%%%%%%%%%%%%%%%%%%%%%%%%%%%%%%%%%%%%%%%%%

We denote by $\eta^+$ (respectively $\eta^-$) the positive
(respectively negative) part of $\eta$.
%%%%%%%%%%%%%%%%%%%%%%%%%%%%%%%%%%%%%%%%%%%%%%%%%%%%%%%%%%%%%%%%%%%%%%%%%%
                 %T H E O R E M  3.2
%%%%%%%%%%%%%%%%%%%%%%%%%%%%%%%%%%%%%%%%%%%%%%%%%%%%%%%%%%%%%%%%%%%%%%%%%%
\begin{Th}
Suppose that for $\theta \in \mathbb{R}^m$ there exists a real
valued nonnegative function $ V_\theta(u)  :  \mathbb{R}^m
\longrightarrow \mathbb{R} $ having  continuous and bounded
partial second derivatives and
\begin{description}
\item[(G1)] $ V_\theta(0)=0,$ and for each  $\ve\in (0, 1),$
$$
\inf_{ \|u\| \geq \ve} V_\theta(u)> 0 ;
$$
\item[(G2)] there exists a set $A\in \mathcal{F}$  with $P^\theta (A)>0$ such that
\; for each $\ve\in (0, 1),$
$$
\sum_{t=1}^\infty \inf_{\ve \le V_\theta(u) \le {1/\ve}}
\left[{\cal N}_t(u)\right]^-=\infty
$$
\end{description}
on $A$, where
\begin{eqnarray}
{\cal N}_t(u)&= &\dot V_\theta(u) \Gm_t^{-1}(\theta+u)b_t(\theta,u)\nonumber\\
&&
+\frac12\sup_{v} \|{\mbox{\"{V}}_\theta}(v)\| E_\theta
\left\{\|\Gm_t^{-1}(\theta+u) \p_t(\theta+u)\|^2 \mid
{\cf}_{t-1}\right\}, \nonumber
\end{eqnarray}
\;\;\;\;\;\;\;
%  u\in {{\mathbb{R}}}^m;
\begin{description}
\item[(G3)]  $\;$ for $\Dl_t=\hat\theta_t-\theta,$
$$
\sum_{t=1}^\infty (1+V_\theta(\Dl_{t-1}))^{-1} \left[{\cal
N}_t(\Dl_{t-1})\right]^+ < \infty, \qquad
 P^\theta\mbox{-a.s.}.
$$
\end{description}

Then $\hat \theta_t \to \theta \;\; (P^\theta$-a.s.)
 for any  initial value
$\hat \theta_0$.
\end{Th}
%%%%%%%%%%%%%%%%%%%%%%%%%%%%%%%%%%%%%%%%%%%%%%%%%%%%%%%%%%%%%%%%%%%%%%%%%%
                    % P R O O F OF T H E O R E M  3.2
%%%%%%%%%%%%%%%%%%%%%%%%%%%%%%%%%%%%%%%%%%%%%%%%%%%%%%%%%%%%%%%%%%%%%%%%%%
 {\bf Proof.} As always (see the convention in Section 2), convergence and all relations between random
variables are meant with probability one w.r.t. the measure
$P^\theta$ unless specified otherwise.
  Rewrite \eqref{rec2}
in the form
$$
\Dl_t=\Dl_{t-1}+ \Gm_t^{-1}(\theta+\Dl_{t-1})
\p_t(\theta+\Dl_{t-1}).
$$
 By the Taylor expansion,
                           %(liap1)
\begin{eqnarray}
V_\theta(\Dl_t) =V_\theta(\Dl_{t-1})+ \dot
V_\theta(\Dl_{t-1})\Gm_t^{-1}(\theta+\Dl_{t-1})
\p_t(\theta+\Dl_{t-1})  \nonumber\\
 +\frac 12 \left[\Gm_t^{-1}(\theta+\Dl_{t-1})
\p_t(\theta+\Dl_{t-1})\right]^T {\mbox{\"{V}}}_\theta(\tilde
\Dl_{t}) \Gm_t^{-1}(\theta+\Dl_{t-1}) \p_t(\theta+\Dl_{t-1}),
\nonumber
\end{eqnarray}
where $\tilde \Dl_{t}\in \mathbb{R}^m$. Taking the conditional
expectation w.r.t. $ {\cf}_{t-1}$ yields
$$
E_\theta\left\{V_\theta(\Dl_t)  \mid {\cal{F}}_{t-1}\right\} \le
V_\theta(\Dl_{t-1})+ {\cal N}_t(\Dl_{t-1}).
$$
Using the obvious decomposition $ {\cal N}_t(\Dl_{t-1})= {[{\cal
N}_t(\Dl_{t-1})]}^+ - {[{\cal N}_t(\Dl_{t-1})]}^-, $ the previous
inequality can be rewritten as
\begin{equation}\label{V}
E_\theta\left\{V_\theta(\Dl_t)  \mid {\cal{F}}_{t-1}\right\} \le
V_\theta(\Dl_{t-1})(1+B_t)+B_t- [{\cal N}_t(\Dl_{t-1})]^-,
\end{equation}
where
$$ B_t=\left(1+V_\theta(\Dl_{t-1}) \right)^{-1}[{\cal
N}_t(\Dl_{t-1})]^+.
$$
 By condition (G3),
\begin{equation}\label{B}
\sum_{t=1}^\infty B_t< \infty.
\end{equation}
According to Lemma A1  in Appendix A (with $X_n=V_\theta(\Dl_n)$,
$\beta_{n-1}=\xi_{n-1}=B_n$ and $\zeta_{n-1}={[{\cal
N}_n(\Dl_{n-1})]}^-$), inequalities \eqref{V} and \eqref{B}  imply
 that the processes $V_\theta(\Dl_{t})$ and
$$
Y_t = \sum_{s=1}^{t} [{\cal N}_s(\Dl_{s-1})]^-
$$
converge to some finite limits. It therefore follows that  $
V_\theta(\Dl_{t}) \to r\ge 0.$ Suppose that $\{r>0\}.$ Then there
exists $\ve > 0$ such that $\ve \le V_\theta(\Dl_{t}) \le {1/\ve}
$ eventually. Because of (G2), this implies that for some
(possibly random) $t_0,$
$$
\sum_{s=t_0}^{\infty} [{\cal N}_s(\Dl_{s-1})]^- \ge
\sum_{s=t_0}^{\infty} \inf_{\ve \le V_\theta(u) \le {1/\ve}}
\left[{\cal N}_s(u)\right]^-=\infty
$$
on the set $A$ with $P^\theta (A)>0,$ which contradicts the existence  of a finite limit of $Y_t.$
Hence, $r=0$   and so,
 $ V_\theta(\Dl_{t}) \to  0$. Now, $\Dl_{t} \to  0$ follows from (G1)
 (otherwise there  would exist a sequence  $t_k \to \infty$  such that
 $\|\Dl_{t_k}\| \ge \ve $ for some $\ve>0,$ and (G1) would imply that
  $ \inf_k V_\theta(\Dl_{t_k}) >  0$).
$ \diamondsuit $

\vskip+0.2cm
%%%%%%%%%%%%%%%%%%%%%%%%%%%%%%%%%%%%%%%%%%%%%%%%%%%%%%%%%%%%%%%%%%%%%%%%%%
%
\noindent {\bf Proof of Theorem 3.1.} As always (see the
convention in Section 2), convergence and all relations between
random variables are meant with probability one w.r.t. the measure
$P^\theta$ unless specified otherwise. Let us show that the
conditions of Theorem 3.1 imply those in Theorem 3.2 with
$V_\theta(u) =(u,u)=u^Tu=\| u\|^2.$ Condition (G1) trivially
holds. Since $\dot V_\theta(u)=2u^T$ and  \"{V}$_\theta(u)=2\times
{\bf 1},$ it follows that
                         % Nt
\begin{equation}\label{Nt}
{\cal N}_t(u)= 2u^T\Gm_t^{-1}(\theta+u)b_t(\theta,u)+ E_\theta
\left\{\|\Gm_t^{-1}(\theta+u) \p_t(\theta+u)\|^2 \mid
{\cf}_{t-1}\right\}.
\end{equation}
Then, by (C1) and (C3),
\begin{eqnarray}\label{Conv}
& & \sum_{t=1}^\infty (1+\|\Dl_{t-1}\|^2)^{-1}
\left[{\cal N}_t(\Dl_{t-1})\right]^+ \nonumber\\
& & \le
  \sum_{t=1}^\infty (1+\|\Dl_{t-1}\|^2)^{-1}
E_\theta \left\{\|\Gm_t^{-1}(\theta+\Dl_{t-1})
\p_t(\theta+\Dl_{t-1}\|^2 \mid {\cf}_{t-1}\right\} \nonumber\\
& & \le
\sum_{t=1}^\infty  B_t < \infty.
\end{eqnarray}
So, (G3) holds. To derive (G2),  using the obvious inequality
$[a]^- \ge -a$   and  (C1), we write
\begin{eqnarray*}
\inf\left[{\cal N}_t(u)\right]^- &\ge &\inf \left[
-2u^T \Gm_t^{-1}(\theta+u)  b_t(\theta,u) \right. \\
&& -E_\theta
\left. \left\{\|\Gm_t^{-1}(\theta+u) \p_t(\theta+u)\|^2 \mid
{\cf}_{t-1}\right\}
\right]\\ \nonumber
&\ge &\inf \left| 2u^T\Gm_t^{-1}(\theta+u)b_t(\theta,u)\right| \\
&&-\sup\left[E_\theta \left\{\|\Gm_t^{-1}(\theta+u)
\p_t(\theta+u)\|^2 \mid {\cf}_{t-1}\right\} \right],
\end{eqnarray*}
where  $\inf$'s and $\sup$'s are taken over $\{u:\ve \le \|u\|^2
\le {1/\ve}\}$.
 From (C3),
$$
\sup\left[E_\theta \left\{\|\Gm_t^{-1}(\theta+u)
\p_t(\theta+u)\|^2 \mid {\cf}_{t-1}\right\} \right] \le
B_t(1+1/{\ve^2})
$$
and $\sum_{t=1}^\infty B_t < \infty.$ Now, using (C2), we finally
obtain
$$
\sum_{t=1}^\infty \inf\left[{\cal N}_t(u)\right]^- \ge
\sum_{t=1}^\infty \inf \left|
2u^T\Gm_t^{-1}(\theta+u)b_t(\theta,u)\right| -
(1+1/{\ve^2})\sum_{t=1}^\infty  B_t=\infty,
 $$
which implies (G2). So, Theorem 3.1 follows on application of
Theorem 3.2. $ \diamondsuit $

\bigskip
%%%%%%%%%%%%%%%%%%%%%%%%%%%%%%%%%%%%%%%%%%%%%%%%%%%%%%%%%%%%%%%%%%%%%%%%%%
%%%%%%%%%%%%%%%%%%%%%%%%%%%%%%%%%%%%%%%%%%%%%%%%%%%%%%%%%%%%%%%%%%%%%%%%%%
                  % R E M A R K 3.3
%%%%%%%%%%%%%%%%%%%%%%%%%%%%%%%%%%%%%%%%%%%%%%%%%%%%%%%%%%%%%%%%%%%%%%%%%%
\begin{rem}
 {\rm It follows from the proof of Theorem 3.2
that if conditions (G1) and (G3) are satisfied then
$(\hat \theta_t - \theta)^2$  converges  ($P^\theta$-a.s.) to a finite limit,
 for any  initial value
$\hat \theta_0$. In particular,  to guarantee this convergence, it suffices
 to require conditions (C1) and (C3) of Theorem 3.1 (this
  can be seen by taking $V_\theta(u) =(u,u)=u^Tu=\| u\|^2$
  and \eqref{Conv}).
}
\end{rem}

%\newpage
                 % 4  SPECIAL MODELS   AND EXAMPLES

\section{SPECIAL MODELS
                   AND EXAMPLES}
\setcounter{equation}{0}
\subsection{The i.i.d. scheme.}
                      %  E X A M P L E 1
 Consider  the classical scheme of
  i.i.d. observations $X_1,X_2,\ldots ,$ with a common probability
density/mass  function
$f(\theta,x), \;\; \theta \in {\mathbb{R}}^m.$
Suppose that $\p(\theta, z)$ is an estimating  function with
$$
\int
\p(\theta,z)f(\theta,z)\mu (dz)=0.
$$
Let us define
the recursive estimator $\hat \theta_t$  by
                             % mleg
\begin{equation}\label{mleg}
\hat \theta_t = \hat \theta_{t-1}  + \frac 1 t\gm^{-1}
(\hat \theta_{t-1})
  \p(\hat \theta_{t-1},X_t),\qquad t\ge 1,
\end{equation}
where $\gm(\theta)$ is a non-random matrix such that
$\gm^{-1}(\theta)$ exists  for any $\theta\in {\mathbb{R}}^m$ and
$\hat\theta_0\in {\mathbb{R}}^m$ is any initial value.

                % C O R O L L A R Y   4.1
\begin{cor}
Suppose that  for any $\theta \in {\mathbb{R}}^m$, the following conditions hold.
\begin{description}
\item[(I)] For any $0<\ve<1,$
$$
\sup_{\ve\le \|u\|\le \frac 1{\ve}}
u^T \;\;  \gm^{-1}(\theta+u)\int \p(\theta+u,x)f(\theta,x)\mu(\,dx)  <  0.
$$
\item[(II)] $\;$For each $ u \in {{\mathbb{R}}^m},$
$$
\int \left \|  \gm^{-1}(\theta+u) \p(\theta+u,x)\right \|^2
f(\theta,x)\mu(\,dx) \leq  K_\theta
(1+\|u\|^2)
$$
for some constant $K_\theta.$
\end{description}
Then  the   estimator
$\hat  \theta_t$   is strongly
consistent for any initial value $\hat\theta_0$.
\end{cor}
\bigskip
\noindent
{\bf Proof} ~~Since
$
   b_t(\theta, u)=b(\theta, u)=\int
\p(\theta+u,z)f(\theta,z)\mu(\,dz)
$
and
$
\Gm_t(\theta)= t\gm(\theta),
$
it is easy to see that  (I) and
 (II) imply   (C1), (C2) and
 (C3) from Theorem 3.1 which yields $(\hat  \theta_t-\theta) \to
 0$ ($P^\theta$-a.s.).

%\vskip+0.2cm
\vskip+0.5cm
Similar  results (for i.i.d. schemes) were obtained
 by  Khas'minskii and Nevelson  \cite{Khas}
 %(when $\p(\theta,x)=l(\theta,x)$ and
 %$\gm(\theta)=i(\theta)$,
Ch.8, $\S$4,  and Fabian \cite{Fab}.
Note  that conditions (I) and (II) are derived from Theorem 3.1
and are sufficient conditions
for the convergence of \eqref{mleg}. Applying
Theorem 3.2 to \eqref{mleg}, one can obtain various
alternative sufficient conditions analogous to those given in
 Fabian (1978). Note also that, in (4.1), the normalising
sequence is $
\Gm_t(\theta)= t\gm(\theta),
$
but Theorems 3.1 and 4.1 allow to consider procedures with arbitrary predictable
$
\Gm_t(\theta).
$
\vskip+0.5cm
%%%%%%%%%%%%%%%%%%%%%%%%%%%%%%%%%%%%%%%%%%%%%%%%%%%%%%%%%%%%%
%       Linear recursion
%%%%%%%%%%%%%%%%%%%%%%%%%%%%%%%%%%%%%%%%%%%%%%%%%%%%%%%%%%%%%
                     %  E X A M P L E 2
\subsection{Linear  procedures.}
Consider the recursion
%          LinRecg
\begin{equation}\label{LinRecg}
\hat\theta_t=\hat\theta_{t-1}+\Gamma_t^{-1}\left( h_t-
\gamma_t\hat\theta_{n-1}\right), \;\;\; t\ge 1,
\end{equation}
where the $\Gamma_t$ and  $\gamma_t$ are predictable processes,
$h_t$ is an adapted  process (i.e., $h_t$ is $\mathcal{F}_t$-measurable for $t\ge 1$)
 and all three are independent of
$\theta.$
The following result  gives a sets of sufficient conditions
for the convergence of \eqref{LinRecg} in the case when the
linear $\psi_t(\theta)=h_t-\gamma_t\theta$ is a martingale-difference.
%$\bf \psi \in \Psi^M.$
%                % C O R O L L A R Y   4.2
\begin{cor}
Suppose that  for any $\theta \in {\mathbb{R}},$
\begin{description}
\item[(a)] $E_\theta\left\{h_t\mid {\cal{F}}_{t-1}\right\}
=\gamma_t\theta,$ ~ for ~$t \ge 1, $~~~~~~~~~~~~~~$P^\theta$-a.s.,
\item[(b)] $0\le \gamma_t/\Gamma_t \le 2-\delta$ eventually for
some $\delta>0, $ ~ and
$$
\sum_{t=1}^\infty {\gamma_t}/{\Gamma_t} =\infty,
$$
on a set $A$  of positive probability $P^\theta$.
\item[(c)]
$$
\sum_{t=1}^\infty \frac
{E_\theta\left\{(h_t-\theta\gamma_t)^2 \mid {\cal{F}}_{t-1}\right\}}
{\Gamma_t^2} < \infty, \qquad P^\theta\mbox{-a.s.}.
$$
\end{description}
 Then   $\hat \theta_t \to \theta \;\;\;( P^\theta$-a.s.) for any  initial value
$\hat \theta_0\in {\mathbb{R}}$ .
\end{cor}
\noindent
{\bf Proof.}
We need to  check that the conditions of Theorem 3.2 hold for for
$V_{\theta}(u)=u^2$.
Using (a) we obtain
$$
b_t(\theta,u)=E_\theta\left\{(h_t-(\theta+u)\gamma_t)\mid
{\cal{F}}_{t-1}\right\}=-u\gamma_t
$$
and
$$
E_\theta\left\{(\psi_t(\theta+u))^2\mid
{\cal{F}}_{t-1}\right\}
=E_\theta\left\{(h_t-(\theta+u)\gamma_t)^2\mid
{\cal{F}}_{t-1}\right\}
$$
$$
=E_\theta\left\{(h_t-\theta\gamma_t)^2\mid
{\cal{F}}_{t-1}\right\}+u^2 \gamma_t^2=\mathcal{P}_t^\theta+u^2 \gamma_t^2,
$$
where $\mathcal{P}_t^\theta=E_\theta\left\{(h_t-\theta\gamma_t)^2\mid
{\cal{F}}_{t-1}\right\}.$
Now, using \eqref{Nt},
$$
N_t(u)=-2 u^2\gamma_t\Gamma_t^{-1}
+\Gamma_t^{-2}\mathcal{P}_t^\theta+u^2\gamma^2_t\Gamma_t^{-2}
$$
$$
=-\delta u^2\gamma_t\Gamma_t^{-1}-u^2\gamma_t\Gamma_t^{-1}\left(
(2-\delta)- \gamma_t\Gamma_t^{-1}\right)
+\Gamma_t^{-2}\mathcal{P}_t^\theta.
$$
To derive (G2),  we use  the obvious inequality
$[a]^- \ge -a$  (for any $a$), conditions  (b) and (c), and write
$$
\sum_{t=1}^\infty
\inf_{\ve \le u^2 \le {1/\ve}}
\left[{\cal N}_t(u)\right]^-\ge
\sum_{t=1}^\infty
\inf_{\ve \le u^2 \le {1/\ve}}
\left(\delta u^2\gamma_t\Gamma_t^{-1} -\Gamma_t^{-2}\mathcal{P}_t^\theta\right)=\infty
$$
on $A$.
To check (G3) we write
$$
\sum_{t=1}^\infty (1+\Dl^2_{t-1})^{-1} \left[{\cal
N}_t(\Dl_{t-1})\right]^+ \le \sum_{t=1}^\infty \left[{\cal
N}_t(\Dl_{t-1})\right]^+\le \sum_{t=1}^\infty
\Gamma_t^{-2}\mathcal{P}_t^\theta<\infty
$$
($P^\theta$-a.s.),  which completes the proof. $ \diamondsuit $

%\appendix
%\section{Appendix}

%\vskip+1.0cm
%\newpage

\bigskip

%\vskip+0.5cm
%%%%%%%%%%%%%%%%%%%%%%%%%%%%%%%%%%%%%%%%%%%%%%%%%%%%%%%%%%%%%%%%%%%%%%%%%%
                  % R E M A R K 4.1
%%%%%%%%%%%%%%%%%%%%%%%%%%%%%%%%%%%%%%%%%%%%%%%%%%%%%%%%%%%%%%%%%%%%%%%%%%
\begin{rem}
{\rm   Suppose that $\Delta\Gamma_t=\gamma_t.$ Then
%         Solg
  \begin{equation}\label{Solg}
\hat\theta_t=\Gamma_t^{-1}\left(\hat\theta_0+\sum_{s=1}^t h_s(X_s)\right).
\end{equation}
This can be easily seen by inspecting the difference $\hat\theta_t-\hat\theta_{t-1}$
for the sequence \eqref{Solg}, to check that \eqref{LinRecg} holds.
It is also interesting to observe that
  since in this case,
 $\Gamma_t=\sum_{s=1}^t \gamma_s,$
$$
\hat\theta_t=\Gamma_t^{-1}\hat\theta_0+\Gamma_t^{-1}\sum_{s=1}^t
\left(h_s(X_s)- \gamma_s\theta\right) +\theta
=\Gamma_t^{-1}\hat\theta_0+\Gamma_t^{-1}M_t^\theta +\theta
$$
where,
$M_t^\theta=\sum_{s=1}^t
\left(h_s(X_s)- \gamma_s\theta\right)$ is a $P^\theta$ martingale. Now, if $\Gamma_t \to \infty$,
 a necessary and sufficient condition for the convergence to $\theta$ is the convergence
 to zero of the sequence $\Gamma_t^{-1}M_t^\theta.$ Condition (c) in Corollary 4.2 is a
 standard sufficient condition in martingale theory to guarantee $\Gamma_t^{-1}M_t^\theta \to 0$  (see
 e.g., Shiryayev  \cite{Shir}, Ch.VII, $\S 5$ Theorem 4).
    The first part of (b) will  trivially hold if $\gamma_t=\Delta \Gamma_t\ge 0$.
Also, in this case,   $\Gamma_t \to \infty$ implies  $
\sum_{t=1}^\infty {\Delta\Gamma_t}/{\Gamma_t} =\infty
$
(see Proposition A3 in Appendix A).
}
\end{rem}
\vskip+0.5cm
%%%%%%%%%%%%%%%%%%%%%%%%%%%%%%%%%%%%%%%%%%%%%%%%%%%%%%%%%%%%%%%%%%%%%%%%%%
                  % R E M A R K 4.2
%%%%%%%%%%%%%%%%%%%%%%%%%%%%%%%%%%%%%%%%%%%%%%%%%%%%%%%%%%%%%%%%%%%%%%%%%%
\begin{rem}
{\rm As a particular example, consider the process
$$
X_t=\theta X_{t-1}+\xi_t, ~~~~~~~~~~~~~~~~~t\ge 1,
$$
where, $\xi_t$ is a $P^\theta$ martingale-difference with
$D_t=E_\theta\left\{\xi_t^2 \mid {\cal{F}}_{t-1}\right\}>0.$
The choice   $h_t=D_t^{-1}X_{t-1}X_t $ and $\Delta\Gamma_t=\gamma_t=D_t^{-1}X_{t-1}^2$,
in \eqref{LinRecg}
 yields the least square estimator of $\theta.$ It is easy to verify that (a) holds. Also, since
$$
E_\theta\left\{(h_t-\gamma_t\theta)^2 \mid {\cal{F}}_{t-1}\right\}
=D_t^{-2}X_{t-1}^2 E_\theta\left\{\xi_t^2 \mid {\cal{F}}_{t-1}\right\}=D_t^{-1}X_{t-1}^2
=\Delta\Gamma_t,
$$
it follows that  (c) in Corollary 4.2 is equivalent to
$
\sum_{t=1}^\infty {\Delta\Gamma_t}/{\Gamma_t^2} <\infty.
$
This, as well as (b) hold if $\Gamma_t \to \infty$
(see Proposition A3 in Appendix A).
So, if $\Gamma_t \to \infty$ the least square procedure is strongly consistent.
If,  e.g., $\xi_t$ are i.i.d. r.v.'s, then $\Gamma_t \to \infty$ for all values of $\theta \in \mathbb{R}$
 (see, e.g, Shiryayev \cite{Shir}, Ch.VII, 5.5).
 }
\end{rem}

                      %  E X A M P L E 3

 \subsection{AR(m) process}
Consider an  AR(m) process
$$
X_i=\theta_{1} X_{i-1}+\dots+\theta_{m} X_{i-m}+\xi_i
=\theta^T X_{i-m}^{i-1}+\xi_i,
$$
where $X_{i-m}^{i-1}=(X_{i-1},\dots,X_{i-m})^T,$
$\theta=(\theta_1,\dots, \theta_m)^T$
 and
${\xi}_i$ is a
sequence of i.i.d. random variables.

A reasonable class of procedures in this model
should have a form
                          %(Arg)
\begin{equation}\label{Arg}
\hat \theta_t=\hat \theta_{t-1}+\Gm_t^{-1}(\hat\theta_{t-1})
 \p_t(X_t-\hat \theta_{t-1}^TX_{t-m}^{t-1}),
\end{equation}
where
$\p_t(z)$ and $\Gm_t^{-1}(z)$ ($z\in {\mathbb{R}}^m$) are respectively  vector and
matrix processes meeting conditions of the previous section.
  Suppose  that  the probability
density function of   ${\xi}_t$ w.r.t. Lebesgue's measure is
$g(x)$. Then the conditional probability density function is
$f_t(\theta, x_t \mid x_1^{t-1})=g(x_t- \theta^T x_{t-m}^{t-1}).$
So, denoting
                          %(Argmle)
\begin{equation}\label{Argmle}
\p_t(z)=-\frac{g'(z)}{g(z)}X_{t-m}^{t-1},
\end{equation}
it is easy to see that
$$
 \p_t(X_t- \theta^T X_{t-m}^{t-1})=
\frac{\dot f_t^T(\theta, X_t \mid X_1^{t-1})}
 {f_ t(\theta, X_t \mid
X_1^{t-1})}
 $$
    and
\eqref{Arg} becomes a likelihood recursive procedure. A possible choice
of  $\Gm_t(z)$  in this case would be
  the conditional Fisher information matrix
 $$
 I_t={i^g}
\sum_{s=1}^t X_{t-m}^{t-1}(X_{t-m}^{t-1})^T
$$
where
$$
i^g=\int  \left(\frac{g'(z)}{g(z)}\right)^2
 g(z) \,dz.
$$
An interesting class of recursive estimators for strongly stationary AR(m)
processes is studied in Campbell \cite{Cam}. These estimators are
recursive versions of robust modifications of the least squares method
and are defined as
%%%%%%%%%%%%%%%%%%%%%%%%%%%%%%%%%%%%%%%%%%%%%%%%%%%%%%%%%%%%%%%%%%%
                          %(Camb)
%%%%%%%%%%%%%%%%%%%%%%%%%%%%%%%%%%%%%%%%%%%%%%%%%%%%%%%%%%%%%
\begin{equation}\label{Camb}
\hat \theta_t=\hat \theta_{t-1}+a_t\gm(X_{t-m}^{t-1})
 \phi(X_t-\hat \theta_{t-1}^TX_{t-m}^{t-1}),
\end{equation}
where $a_t$ is a sequence of a positive numbers with $a_t\to 0,$
$\phi$ is a bounded scalar function and
$\gm(u)$ is a vector function of the form $u h(u)$ for some
non-negative function $h$ of $u$ (See also Leonov \cite{Leon}).
The class of procedures of type \eqref{Camb} is clearly  a subclass
of that defined by \eqref{Arg}  and therefore can be studies
 using the results of the previous section.

 Suppose that $\xi_i$ are i.i.d.
random variables with  a bell-shaped,  symmetric
about  zero probability density function $g(z)$ (that is,
$g(-z)=g(z),$ and $g\downarrow 0$ on ${\mathbb{R}}_+$).
  Suppose also that $\phi(x)$ is
an odd,  continuous in zero function.
Let  us write  conditions of Theorem 3.1  for
%                         (psi)
\begin{equation}\label{psi}
\Gamma(\theta)=a_t^{-1}{\bf 1} ~~~\mbox{and}~~~\p_t(\theta)=
X_{t-m}^{t-1}h\left(X_{t-m}^{t-1}\right)
\phi\left(X_t-\theta^TX_{t-m}^{t-1}\right).
\end{equation}
We have
\begin{eqnarray}
E_\theta\left\{
\phi\left(X_t-(\theta+u)^TX_{t-m}^{t-1}\right)
\mid {\cal{F}}_{s-1}
\right\}
& & =
E_\theta\left\{
\phi\left(\xi_t-u^TX_{t-m}^{t-1}\right)
\mid {\cal{F}}_{s-1}
\right\}
\nonumber\\
& & =\int
\phi\left(z-u^TX_{t-m}^{t-1}\right)
g(z)dz. \nonumber
\end{eqnarray}
It follows from  Lemma A2 in  Appendix A that if
$w\not =0,$
$$
G(w)=-w\int_{-\infty}^{\infty}\phi\left(z-w\right) g(z)dz>0.
$$
Therefore,
 \begin{eqnarray}\label{Negat}
u^T\Gm_t^{-1}(\theta +u)b_t(\theta,u)
=
a_t u^T X_{t-m}^{t-1}h(X_{t-m}^{t-1})E_\theta\left\{
\phi(\xi_t-u^TX_{t-m}^{t-1})
\mid {\cal{F}}_{s-1}
\right\}
\nonumber
\end{eqnarray}
\begin{eqnarray}
=-a_t~ h\left(X_{t-m}^{t-1}\right)G(u^TX_{t-m}^{t-1})\le 0.
\end{eqnarray}
Also,  since $\phi$ is a bounded function,
$$
E_\theta \left\{\|\Gm_t^{-1}(\theta+u)\p_t(\theta+u)\|^2\mid {\cf}_{t-1}\right\}
 \le C^\theta a_t^2 \|X_{t-m}^{t-1}\|^2h^{2} (X_{t-m}^{t-1})
 $$
for some positive constant $C^\theta$. Therefore, conditions
of Theorem 3.1 hold if ($P^\theta$-a.s.),
%                         (ArC1)
\begin{equation}\label{ArC1}
\sum_{t=1}^\infty ~a_t h\left(X_{t-m}^{t-1}\right)\inf_{\ve \le \|u\| \le
{1/\ve}}
 G(u^TX_{t-m}^{t-1})=\infty
 \end{equation}
 and
                            %(ArC2)
 \begin{equation}\label{ArC2}
\sum_{t=1}^\infty a_t^2 \|X_{t-m}^{t-1}\|^2h^{2} (X_{t-m}^{t-1})<\infty.
\end{equation}
If $X_t$ is a  stationary process,  these conditions can be verified
using limit theorems for  stationary processes.
Suppose, e.g., that
 $a_t=1/t$,  $~~h({\bf x})\not= 0$ for any  ${\bf x}\not=0$, and $g(z)$ is continuous.  Then
$~~h({\bf x})\inf_{\ve \le \|u\| \le {1/\ve}} G(u^T{\bf x})>0$ for any ${\bf x}\neq
0$ (see Appendix A, Lemma A2). Therefore,  it follows from an ergodic theorem for
stationary processes that in probability $P^\theta,$
                           %(Stat)
 \begin{equation}\label{Stat}
\lim_{t\to \infty}\frac1t\sum_{s=1}^t h\left(X_{s-m}^{s-1}\right)
\inf_{\ve \le \|u\| \le {1/\ve}}
G(u^TX_{s-m}^{s-1}) ~~~~ >0.
\end{equation}

Now, \eqref{ArC1} follows from
Proposition A4, in Appendix A.

Examples of the procedures of type \eqref{Camb} as well as
some simulation results are presented in Campbell \cite{Cam}.

\bigskip

%%%%%%%%%%%%%%%%%%%%%%%%%%%%%%%%%%%%%%%%%%%%%%%%%%%%%%%%%%%%%%%%%%%%%%%%%%
                  % Example 4.2
%%%%%%%%%%%%%%%%%%%%%%%%%%%%%%%%%%%%%%%%%%%%%%%%%%%%%%%%%%%%%%%%%%%%%%%%%%
\subsection{An explicit  example} As a particular example of (4.4), consider the process
$$
X_t=\theta X_{t-1}+\xi_t, ~~~~~~~~~~~~~~~~~t\ge 1,
$$
where, $\xi_t,~~~ t\ge 1,$ are independent Student random variables with degrees of
freedom $\alpha$. So, the probability density functions of $\xi_t$ is
$$
g(x)=C_\alpha \left(
1+\frac {x^2} {\alpha}
\right)^{-\frac{\alpha+1}2}
$$
where $C_\alpha={\bf\Gamma}((\alpha+1)/2)/(\sqrt{\pi\alpha} ~{\bf\Gamma}(\alpha/2)).$

Since
$$
\frac{g'(z)}{g(z)}=
-(\alpha+1)\frac{z}{\alpha+z^2}
$$
(see also \eqref{Argmle}),
$$
\frac{\dot f_t(\theta, X_t \mid X_{t-1})}
 {f_ t(\theta, X_t \mid
X_{t-1})}
 =-X_{t-1}\frac{g'}{g}(X_t- \theta X_{t-1})=
 (\alpha+1)X_{t-1}\frac{X_t- \theta X_{t-1}}{\alpha+({X_t- \theta X_{t-1}})^2}
$$
    and
 the conditional Fisher information is
 $$
 I_t={i^g}
\sum_{s=1}^t X_{t-1}^2
$$
where
\begin{eqnarray}
i^g
=\int  \left(\frac{g'(z)}{g(z)}\right)^2
 g(z) \,dz & &=C_\alpha (\alpha+1)^2
\int  \frac{z^2\,dz}{(\alpha+z^2)^2(1+\frac{z^2}{\alpha})^{\frac{\alpha+1}{2}}} \nonumber
\\
& &=
C_\alpha \frac{(\alpha+1)^2}{\sqrt{\alpha}}
\int  \frac{z^2\,dz}{(1+z^2)^{\frac{\alpha+5}{2}}}\nonumber
\\
& &
=C_\alpha \frac{(\alpha+1)^2}{\sqrt{\alpha}}
\frac{\sqrt{\pi}{\bf\Gamma}((\alpha+5)/2-3/2)}{2{\bf\Gamma}((\alpha+5)/2)}
\nonumber
\\
& & =\frac{2(\alpha+1)}{\alpha+3}.
\nonumber
\end{eqnarray}
Therefore, a likelihood recursive procedure  is
                         %(Stud)
\begin{equation}\label{Stud}
\hat \theta_t=\hat \theta_{t-1}+I_t^{-1}(\hat\theta_{t-1})
 (\alpha+1)X_{t-1}\frac{X_t- \hat \theta_{t-1}X_{t-1}}{\alpha+({X_t- \hat \theta_{t-1}X_{t-1}})^2},
 ~~~~~t\ge 1,
\end{equation}
where $\hat \theta_0$ is any starting point. Note that $I_t$ can also be
derived recursively by
$$I_t=I_{t-1}+{i^g} X_{t-1}^2.$$
Clearly, {\eqref{Stud}} is a recursive procedure of type {\eqref{Camb}} but with
a stochastic normalizing  sequence $a_t=I_t^{-1}$. Now,  $\psi_t$
is of a form of {\eqref{psi}} with
$h(u)=1$ and $\phi(z)=(\alpha+1){z}/{(\alpha+z^2)},$
and $g(z)$ is a bell-shaped and   symmetric
about  zero. Therefore, to show  convergence to $\theta,$ it suffices
 to check conditions {\eqref{ArC1}} and {\eqref{ArC2}}, which, in this case
 can be written as
 %                         (StArC1)
\begin{equation}\label{StArC1}
\sum_{t=1}^\infty ~\frac{1}{I_t} \inf_{\ve \le |u| \le
{1/\ve}}
 G(u X_{t-1})=\infty
 \end{equation}
 and
                            %(StArC2)
 \begin{equation}\label{StArC2}
\sum_{t=1}^\infty \frac  {X_{t-1}^2}{I_t^2} <\infty,
\end{equation}
($P^\theta$-a.s.). We have,
$I_t \to \infty $ for any $\theta\in \mathbb{R}$
 (see, e.g, Shiryayev \cite{Shir}, Ch.VII, 5.5). Since $\Delta I_t=i^g(X_{t-1})^2,$
 we obtain that {\eqref{StArC2}} follows from Proposition A3 in Appendix A.
Let us assume now that $|\theta| <1.$ By Lemma A2 in Appendix A, $\inf_{\ve \le |u| \le {1/\ve}}
G(ux)>0$ for any $x\not= 0.$ Then  if we assume that the the process is strongly stationary, it follows from
the ergodic theorem that in probability $P^\theta$,
$$
\lim_{t\to\infty} \frac 1t I_t~~ >0
~~~~~ \mbox{and} ~~~~~
\lim_{t\to\infty}\frac 1t \sum_{s=1}^t
\inf_{\ve \le |u| \le {1/\ve}}
G(uX_{s-1})~~~ >0.
$$
(It can be proved that these hold without assumption of strong stationarity.)
Therefore, in probability $P^\theta,$
$
\lim \frac 1{I_t} \sum_{s=1}^t
\inf_{\ve \le |u| \le {1/\ve}}
G(uX_{s-1})~ >0
$
and {\eqref{StArC1}} now follows on application of Proposition A4 in Appendix A.

%%%%%%%%%%%%%%%%%%%%%%%%%%%%%%%%%%%%%%%%%%%%%%%%%%%%%%%%%%%%%%%%%%%%%%%%%
                  % R E M A R K 4.3
%%%%%%%%%%%%%%%%%%%%%%%%%%%%%%%%%%%%%%%%%%%%%%%%%%%%%%%%%%%%%%%%%%%%%%%%%%
\begin{rem}
{\rm We have shown above that the recursive estimator
{\eqref{Stud}} is strongly consistent, i.e., converges to $\theta$ a.s., if $|\theta| < 1.$ It is worth mentioning that
 {\eqref{StArC2}}, and therefore, {\eqref{ArC2}} holds for any  $\theta \in \mathbb{R}$,
which guarantees (C3)  of Theorem 3.1. Also,   {\eqref{Negat}} implies that  (C1) of Theorem 3.1 holds as well.
Therefore, according to  Remark 3.3,  we obtain that $|\hat \theta_t-\theta|$  converges ($P^\theta$-a.s.)
to a finite limit
for any $\theta \in \mathbb{R}$.
}
\end{rem}

\begin{figure}
\begin{center}
\resizebox{0.9\textwidth}{!}{\includegraphics{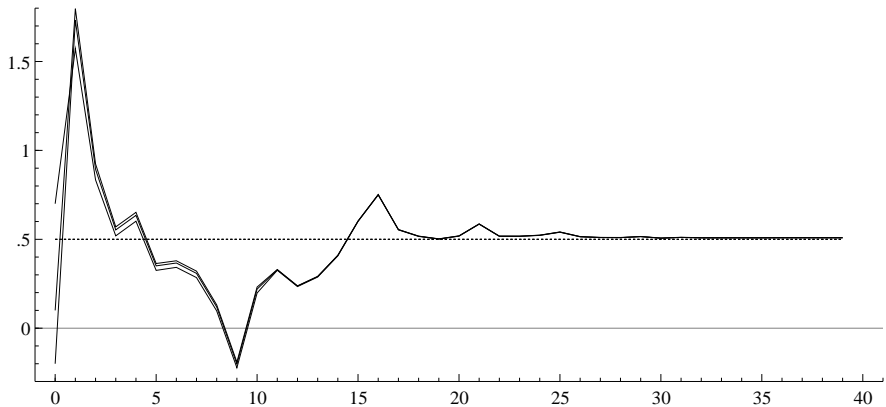}}
\end{center}
\caption{\small Realisations
of {\eqref{Stud}}   for $\alpha=3$ and $\theta=0.5$ for three different
starting values $\theta_0=-0.2, ~ 0,1$ and $0.7$. The number of observations is $40$.}
\end{figure}
%%%%%%%%%%%%%%%%%%%%%%%%%%%%%%%%%%%%%%%%%%%%%%%%%%%%%%%%%%%%%%%%%%%%%%%%%
                  % R E M A R K 4.4
%%%%%%%%%%%%%%%%%%%%%%%%%%%%%%%%%%%%%%%%%%%%%%%%%%%%%%%%%%%%%%%%%%%%%%%%%%
\begin{rem}
{\rm
Note that  conditions {\eqref{StArC1}} and {\eqref{StArC2}} will still hold if we replace
$I_t$ by $c_tI_t$ where $c_t$ is a sequence of non-negative r.v.'s
such that $c_t =1$  eventually.
So, the procedure {\eqref{Stud}} will remain consistent  if $I_t$ is replaced
 by $c_tI_t,$ i.e., if   tuning constants are introduced.
We have shown that the procedure is consistent, i.e., the recursive
estimator is close to the value of the unknown parameter for the large $t$'s. But in practice,
 the tuning constants may be useful to control the behaviour
of a recursion  at the ``beginning'' of the procedure.
 Fig.1 shows    realisations
of  {\eqref{Stud}} for $\alpha=3$ and $\theta=0.5$ for three different
starting values.  The number of observations is 40. As we can see from these graphs, the recursive
procedure, at each step moves in the direction of the parameter (see also Remark 3.2), but
oscillates quite violently for  the first ten steps and then settles down nicely
after another ten steps. This oscillation is due to the small values of the normalising sequence
for the first several steps and
 can be dealt with by introducing tuning constants. On other occasions, it
 may be desirable to lower the value of the normalising sequence for the first several
 steps. This happens when a
 procedure settles down too quickly without any, or little
 oscillation (before reaching the actual value of the parameter).
The detailed discussion of these and related topics will appear elsewhere.
}
\end{rem}

\vskip+1.0cm

%\newpage
\begin{center}
{\large \bf APPENDIX A}
\end{center}

%\appendix
%\section{Appendix}

\noindent
{\bf Lemma  A1}
{\it  Let ${\cal F}_0,{\cal F}_1, \dots$ be a non-decreasing sequence of
$\sigma$-algebras and  $X_n, \beta_n, \xi_n, \zeta_n \in {\cal F}_n,
\;\; n\ge 0,$ are nonnegative r.v.'s such that
$$
E(X_n|{\cal F}_{n-1}) \le X_{n-1}(1+\beta_{n-1})+\xi_{n-1}-
\zeta_{n-1}, \;\;\; n\ge 1
$$
eventually. Then}
$$
\{\sum_{i=1}^{\infty} \xi_{i-1}<\infty\}\cap
\{\sum_{i=1}^{\infty} \beta_{i-1} <\infty\}\subseteq
\{X\rightarrow\}\cap\{
\sum_{i=1}^{\infty} \zeta_{i-1} <\infty\}
\quad (P\mbox{-}a.s.),
$$
where $\{X\rightarrow\}$ denotes the set where $\lim_{n\to \infty} X_n$ exists and is
finite.
\vskip+0.2cm
\noindent
{\bf Remark} Proof can be found in
 Robbins and  Siegmund \cite{Rob2}. Note also that this lemma is a  special case of
the theorem on the convergence sets nonnegative semimartingales
(see, e.g., Lazrieva et al \cite{Laz1}).

\bigskip
\noindent
{\bf Lemma  A2} {\it Suppose that $g \not\equiv 0$
is a nonnegative even  function
on $\mathbb{R}$
and  $g\downarrow 0$ on $\mathbb{R}_+$.
Suppose also that $\phi$ is
a measurable odd function on $\mathbb{R}$ such that $\phi(z) > 0$
for $z>0$
and $\int_{\mathbb{R}}|\phi(z-w)| g(z)dz < \infty$
for all $w \in \mathbb{R}$.  Then
$$
w\int_{-\infty}^{\infty}\phi\left(z-w\right) g(z)dz < 0
\leqno{(A1)}
$$
for any
$w\not =0.$ Furthermore, if $g(z)$ is continuous, then for any} $\ve\in (0, 1)$
$$
\sup_{\ve \le |w| \le 1/\ve} w\int_{-\infty}^{\infty}\phi\left(z-w\right)
g(z)dz < 0 . \leqno{(A2)}
$$
\noindent
{\bf Proof} Denote
$$
\Phi(w)=\int_{-\infty}^{\infty}\phi\left(z-w\right) g(z)dz=\int_{-\infty}^{\infty}\phi(z) g(z+w)dz.
\leqno{(A3)}
$$
Using the change of variable $z \longleftrightarrow -z$ in the
integral over $(-\infty, 0)$ and
the equalities $\phi(-z)=-\phi(z)$ and $g(-z+w) = g(z-w)$,
we obtain
\begin{eqnarray}
\Phi(w)
& & =\int_{-\infty}^{0}\phi(z) g(z+w)dz
+\int_{0}^{\infty}\phi(z) g(z+w)dz
\nonumber \\
& & = \int_{0}^{\infty}\phi(z)\left(g(z+w)-g(-z+w)\right)dz
\nonumber \\
& & = \int_{0}^{\infty}\phi(z)\left(g(z+w)-g(z-w)\right)dz .
\nonumber
\end{eqnarray}
Suppose now that $w>0$. Then $z-w$ is closer to $0$ than
$z+w$, and the properties of $g$ imply that
$g(z+w)-g(z-w)\le 0$.
Since $\phi(z) > 0$ for $z>0$, $\Phi(w) \le 0$.
The equality $\Phi(w) = 0$ would imply that
$g(z+w)-g(z-w)= 0$ for {\it all}  $z \in (0, +\infty)$
since, being monotone, $g$ has right and left
limits at each point of $(0, +\infty)$. The last equality,
however, contradicts the restrictions on $g$. Therefore,
(A1) holds.
Similarly, if  $w<0$, then $z+w$ is closer to $0$ than
$z-w$, and  $g(z+w)-g(z-w) \ge 0$. Hence
$w\left(g(z+w)-g(z-w)\right) \le 0$, which yields (A1)
as before.

To prove (A2) note that the continuity of $g$ implies
that $g(z+w)-g(z-w)$ is a continuous functions of
$w$ and (A2) will follow from (A1) if one proves that $\Phi(w)$ is
also continuous in $w$.  So, it is sufficient to show that
the integral in (A3) is uniformly convergent for
$\ve \le |w| \le 1/\ve$. It follows from the restrictions
we have placed on $g$ that there exists $\delta > 0$ such that
$g \ge \delta$ in a neighbourhood of 0. Then the condition
$$
\int_{0}^{\infty}\phi(z)\left(g(z+w) + g(z-w)\right)dz =
\int_{-\infty}^{\infty} |\phi(z-w)| g(z)dz < \infty , \ \
\forall w \in \mathbb{R}
$$
implies that $\phi$ is locally integrable on $\mathbb{R}$.
It is easy to see that, for any $\ve\in (0, 1)$,
$$
g(z \pm w) \le g(0)\chi_\varepsilon(z) + g(z - 1/\ve) ,
\ \ z \ge 0,  \ \ \ve \le |w| \le 1/\ve ,
$$
where $\chi_\varepsilon$ is the indicator function of the
interval $[0, 1/\ve]$. Since the function
$\phi(\cdot)\,(g(0)\chi_\varepsilon + g(\cdot - 1/\ve))$ is
integrable on $(0, +\infty)$ and does not depend on $w$,
we conclude that
the integral in (A3) is indeed uniformly convergent for
$\ve \le |w| \le 1/\ve$. $\diamondsuit$
%\newpage

\vskip+0.3cm
\noindent
{\bf Proposition A3} {\it
If  $d_n$ is a  nondecreasing sequence  of positive numbers such that $d_n\to +\infty$,
then
$$
\sum_{n=1}^\infty \tr d_n/d_n=+\infty
$$
 and
$$
\sum_{n=1}^\infty \tr d_n/d_n^2 <+\infty.
$$
}
\vskip+0.2cm
\noindent
{\bf Proof} The first claim is easily obtained by contradiction from the Kronecker lemma
(see, e.g., Lemma 2, $\S$3, Ch. IV in Shiryayev \cite{Shir}). The second one is proved by the
following argument
$$
0 \le \sum_{n=1}^N \frac{\tr d_n}{d_n^2} \le
\sum_{n=1}^N \frac{\tr d_n}{d_{n - 1} d_n} =
\sum_{n=1}^N \left(\frac{1}{d_{n - 1}} - \frac{1}{d_n}\right) =
\frac{1}{d_0} - \frac{1}{d_N} \to \frac{1}{d_0} <+\infty.
$$
$\diamondsuit$
%\newpage

\vskip+0.3cm
\noindent
{\bf Proposition A4}
{\it
Suppose that  $d_n, ~ c_n,$ and $c$ are random variables, such that, with probability 1,
$d_n>0, ~ c_n\ge 0, ~ c>0$ and $d_n\to +\infty$ as $n \to \infty.$
Then
$$
\frac 1{d_n}\sum_{i=1}^n  c_i \to c \quad \mbox{in probability}
$$
implies
$$
\sum_{n=1}^\infty \frac {c_n}{d_n}=\infty \quad \mbox{with probability 1}.
$$
}

\vskip+0.3cm
\noindent
{\bf Proof}
 Denote
$\xi_n=\frac 1{d_n}\sum_{i=1}^n  c_i$.
Since $\xi_n \to c$  in probability,
it follows that there exists a
subsequence $\xi_{i_n}$ of  $\xi_n$ with the property that  $\xi_{i_n} \to c$
with probability 1.
Now,  assume that
$
\sum_{n=1}^\infty  {c_n}/{d_n}<\infty
$
on a set $A$ of  positive probability. Then,
 it follows from the Kronecker lemma,
(see, e.g., Lemma 2, $\S$3, Ch. IV in Shiryayev \cite{Shir}) that $
\xi_n \to 0 $ on $A$. Then it follows that $\xi_{i_n} \to 0$ on
$A$ as well, implying that $c=0$ on $A$ which contradicts the
assumptions that $c>0$ with probability 1. $\diamondsuit$

%\end{article}
%
\end{document}